\documentclass{amsart}

\newtheorem{theorem}{Theorem}[section]
\newtheorem{lemma}[theorem]{Lemma}

\newtheorem{proposition}[theorem]{Proposition}

\theoremstyle{definition}
\newtheorem{definition}[theorem]{Definition}
\newtheorem{example}[theorem]{Example}

\theoremstyle{remark}
\newtheorem{remark}[theorem]{Remark}

\numberwithin{equation}{section}

\title{On a special class of smooth codimension two subvarieties in $\Pn$, $n \geq 5$}
\author{C. Folegatti}
\date{08/06/2004\\
The author is supported by MIUR project "Geometry on Algebraic
Varieties"}

\newcommand{\epf}{\ensuremath{\diamondsuit}}
\newcommand{\id}{\ensuremath{\mathcal{I}}}
\newcommand{\oc}{\ensuremath{\mathcal{O}}}
\newcommand{\fc}{\ensuremath{\mathcal{F}}}

\newcommand{\ec}{\ensuremath{\mathcal{E}}}

\newcommand{\nc}{\ensuremath{\mathcal{N}}}

\newcommand{\lc}{\ensuremath{\mathcal{L}}}
\newcommand{\rc}{\ensuremath{\mathcal{R}}}

\newcommand{\bc}{\ensuremath{\mathcal{B}}}
\newcommand{\Pt}{\mathbb{P}^3}

\newcommand{\Pq}{\mathbb{P}^4}
\newcommand{\Pcq}{\mathbb{P}^5}
\newcommand{\Psx}{\mathbb{P}^6}
\newcommand{\Pn}{\mathbb{P}^n}
\newcommand{\Pun}{\mathbb{P}^1}

\newcommand{\bZ}{\mathbb{Z}}

\newcommand{\sig}{\Sigma}
\newcommand{\de}{\Delta}

\begin{document}
\maketitle

\section{Introduction}
We work on an algebraically closed field of characteristic zero.\\
By Lefschetz's theorem, a smooth codimension two subvariety $X \subset \Pn$, $n \geq 4$, which is not a complete intersection, lying on a hypersurface $\sig$, verifies $dim(X \cap Sing(\sig)) \geq n-4$.\\
In this paper we deal with a situation in which the singular locus of $\sig$ is as large as can be, but, at the same time, the simplest possible: we assume $\sig$ is an hypersurface of degree $m$ with an ($m$-2)-uple linear subspace of codimension two.\\
More generally, we are concerned with smooth codimension two subvarieties $X \subset \Pn$, $n \geq 5$.\\
In the first part we consider smooth subcanonical threefolds $X \subset \Pcq$ and we prove that if $deg(X) \leq 25$, then $X$ is a complete intersection (Prop. \ref{prop}). In the second section we study a particular class of codimension two subvarieties and we prove the following result.
\begin{theorem}
\label{th}
Let $X \subset \Pn$, $n \geq 5$, be a smooth codimension two subvariety (if $n=5$ assume $Pic(X)=\bZ H$) lying on a hypersurface $\sig$ of degree $m$, which is singular, with multiplicity $m-2$, along a linear subspace $K$ of dimension $n-2$. Then $X$ is a complete intersection.
\end{theorem}
This gives further evidence to Hartshorne conjecture in codimension two.\\
It is enough to prove the theorem for $n=5$, the result for higher dimensions will follow by hyperplane sections. For $n=5$ it is necessary to suppose $Pic(X)=\bZ H$, whereas for $n \geq 6$, thanks to Barth's theorem, this hypothesis is always verified.\\
The proof for $n=5$ goes as follows. Using the result of the first part we may assume $d \geq 26$, then we prove, under the special assumptions of the theorem, that either $deg(X)$ is less than $25$ or we use the result of Lemma \ref{lem-sub} to conclude that $S$ is a complete intersection.\\
By the way we give a little improvement of earlier results on the non existence of rank two vector bundles on $\Pq$ with small Chern classes, see Lemma \ref{lem-E}.\\
I am very grateful to Ph. Ellia for many useful suggestions and
for his support. I would like also to thank the referee for his remarks, which helped me to improve this work.

\section{Smooth subcanonical threefolds in $\Pcq$}
Let $X$ be a smooth subcanonical threefold in $\Pcq$, of degree $d$, with $\omega_X \cong \oc_X(e)$. Let $S = X \cap H$ be the general hyperplane section of $X$, $S$ is a smooth subcanonical surface in $\Pq$, indeed by adjunction it is easy to see that $\omega_S \cong \oc_S(e+1)$. Again we set $C$ the general hyperplane section of $S$, $C$ is a smooth subcanonical curve in $\Pt$, with $\omega_C \cong \oc_C(e+2)$.\\
We can compute the sectional genus $\pi(S)$, indeed since
$\omega_C \cong \oc_C(e+2)$ it follows that $\pi = g(C)= 1 +
\frac{d(e+2)}{2}$.
\begin{lemma}
\label{q}
With the notations above, $q(S)=0$ and all hyperplane sections $C$ of $S$ are linearly normal in $\Pt$.
\end{lemma}
\textit{Proof:} By Barth's theorem we know that if $X \subset \Pcq$ is a smooth threefold, then $h^1(\oc_X)=0$. Let us consider the exact sequence: $0 \to \oc_X(-1) \to \oc_X \to \oc_S \to 0$. By taking cohomology and observing that $h^2(\oc_X(-1))=h^1(\omega_X(1))=0$ by Kodaira, we get the result.\epf \\
\\
If we look at the surface $S$, we can observe that most of its invariants are known. Hence it seems natural to consider the double points formula in order to get some more information.\\
Since $q(S)=0$, $\pi-1=\frac{d(e+2)}{2}$ and $K^2=d(e+1)^2$, the formula becomes: $d(d-2e^2-9e-17) = -12(1+p_g(S))$, where the quantity $1+p_g(S)$ is strictly positive. We have the following condition:
$$d(d-2e^2-9e-17) \equiv 0 \:\:(mod \:\: 12) \:\:\:\:\:\:\:\:\:\:(1)$$
\begin{proposition}
\label{prop}
Let $X \subset \Pcq$ be a smooth subcanonical threefold of degree $d$, then if $d \leq 25$, $X$ is a complete intersection.
\end{proposition}
\textit{Proof:} We recall that for a smooth subcanonical threefold in $\Pcq$ with $\omega_X \cong \oc_X(e)$ we have $e \geq 3$, unless $X$ is a complete intersection (see \cite{BC}). Let $G(d,3)=1+ \frac{d(d-3)-2r(3-r)}{6}$ be the maximal genus of a curve of $\Pt$ of degree $d=3k+r$, $0 \leq r \leq 2$, not lying on a surface of degree two.
If we compare the value of $\pi$ computed before with this (using $e \geq 3$), we see that if $d \leq 17$, then $h^0(\id_C(2)) \neq 0$. Since by Severi's and Zak's theorems on linear normality $h^1(\id_S(1))=h^1(\id_X(1))=0$, it follows that $h^0(\id_X(2)) \neq 0$ and this implies that $X$ is a complete intersection (see \cite{EF}, Theorem 1.1).\\
If $d=18$, then $\pi=G(18,3)$. It follows that $C$ is a.C.M. then by the exact sequence: $0 \to \oc_{\Pt} \to \ec \to \id_C(e+6) \to 0$ we obtain $h^1(\ec(k))=0$ for all $k \in \bZ$. Hence by Horrocks' theorem $\ec$ is split and then $C$ is a complete intersection. Since this holds for the general $\Pt$ section $C$, the same holds for $S$ and for $X$.\\
If $d=19$ then $C$ lies on a quadric surface unless $\pi=1+\frac{19(e+2)}{2} \leq G(19,3)$. This inequality yields $e=3$ but if we look at formula $(1)$ we see that this is not possible.\\
If $d=20$, then $\pi=G(20,4)$, $C$ is a.C.M. and we argue as in the case $d=18$ to conclude that $X$ is a complete intersection.\\
If $d=21,22,23$ and if $h^0(\id_C(4)) \neq 0$, then thanks to the "lifting theorems" in $\Pq$ and $\Pcq$ (see \cite{Me}) we have $h^0(\id_X(4)) \neq 0$ and again by \cite{EF} $X$ is a complete intersection. We then assume $h^0(\id_C(4))=0$ and using the fact that $\pi = 1 +  \frac{d(e+2)}{2} \leq G(d,5)$, we obtain $e=3$. However this is not possible because of formula $(1)$.\\
If $d=24$ we still get $e=3$, but formula $(1)$ is satisfied. We have the following exact sequence $0 \to \oc \to \ec(5) \to \id_X(9) \to 0$, where $\ec$ is a rank two vector bundle with $c_1(\ec)=-1$ and $c_2(\ec)=4$. If $h^0(\id_X(4))=0$, then $h^0(\ec)=0$, which is not possible since by \cite{D} there exists no rank two stable vector bundle with such Chern classes. Hence it would be $h^0(\id_X(4)) \neq 0$ and this implies (see \cite{EF}, Theorem 1.1) that $X$ is a complete intersection but this is also impossible since the system given by the equations $a+b=-1$ and $ab=4$ does not have solution in $\bZ$.\\
If $d=25$, supposing $h^0(\id_C(4))=0$ we obtain $e=4$. In that case we have exactly $\pi=G(25,5)=76$ and this means that if
$h^0(\id_C(4)) =0$, then $C$ is a.C.M.. It follows that $C$, and then $X$, is a complete intersection. \epf
\begin{remark}
If we perform the same calculations of the proof of \ref{prop} for $d=26$, we have that $e=3$.\\
Now if we consider subcanonical threefolds in $\Pcq$ with $e=3$, by Kodaira we have that $h^0(\oc_X(4))= \chi(\oc_X(4))$. By Riemann-Roch formula for threefolds (see \cite{BC}) we compute $\chi(\oc_X(4))= \frac{5d(50-d)}{24}$. Since $h^0(\oc_{\Pcq}(4))=126$, it is easy to see that for $d \geq 30$ it must be $h^0(\id_X(4)) \neq 0$, hence $X$ is a complete intersection.\\
On the other hand for $26 \leq d \leq 30$ the unique value of $d$ satisfying $(1)$ is $d=26$.
Thus we have shown that, among smooth threefolds in $\Pcq$ with $e=3$, the only possibility for $X$ not to be a c.i. is if $d=26$.
\end{remark}
We conclude this section with some result about rank two vector bundles. Let us start with a lemma concerning subcanonical double structures.
\begin{lemma}
\label{lem-Z}
Let $Y \subset \Pn$, $n \geq 4$, be a complete intersection of codimension two. Let $Z$ be a l.c.i. subcanonical double structure on $Y$. Then if $emdim(Y) \leq n-1$, $Z$ is a complete intersection.
\end{lemma}
\textit{Proof:} By \cite{MM} we have that any doubling of a l.c.i. $Y$ with $emdim(Y) \leq dim(Y)+1$ is obtained by the Ferrand construction. Hence there is a surjection $\nc^{\vee}_Y \to \lc \to 0$ where $\lc$ is a locally free sheaf of rank one on $Y$. Taking into account that $\omega_{Z|Y} \cong \omega_Y \otimes \lc^{\vee}$ (see \cite{BF}) and recalling that $Z$ is subcanonical and $Y$ is a c.i., we obtain that $\lc \cong \oc_Y(l)$ for a certain $l \in \bZ$.\\
On the other hand, since $Y$ is a complete intersection, say
$Y=F_a \cap F_b$, we have $\nc_Y \cong \oc_Y(a) \oplus \oc_Y(b)$,
then the sequence above becomes: $\oc_Y(-a) \oplus \oc_Y(-b)
\stackrel{f}{\to} \oc_Y(l) \to 0$. The map $f$ is given by two
polynomials of degree respectively $a+l$ and $b+l$. If $F$ and $G$
are both not constant, it follows, since $n \geq 4$, that
$B:=(F)_0 \cap (G)_0 \cap Y \neq \emptyset$. For each $x \in B$
the induced map $f_x$ on the stalks is not surjective: absurd.
Thus necessarily $F$ or $G$ is a non zero constant, i.e. either
$l=-a$ or $l=-b$. If $l=-a$ (resp. $l=-b$) we are doubling $Y$ on
$F_b$, $Z=F_a^2 \cap F_b$ (resp. we are doubling $Y$ on $F_a$, $Z=
F_a \cap F_b^2$). In any case, $Z$ is a complete intersection.\epf
\begin{lemma}
\label{lem-Z2}
Let $Z \subset \Pq$ be a l.c.i. quartic surface with $\omega_Z \cong \oc_Z(-a)$. If $a \geq 3$, then $Z$ is a complete intersection.
\end{lemma}
\textit{Proof:} Let $C$ be the hyperplane section of $Z$ and let
$C_{red}=\tilde{C}_1 \cup \ldots \cup \tilde{C}_s$ be the
decomposition of $C_{red}$ in irreducible components, hence $C=C_1
\cup \ldots \cup C_s$, where $C_i$ is a multiple structure on
$\tilde{C}_i$ for all $i$ . We have $\omega_C \cong \oc_C(-a+1)$, on
the other hand $\omega_{C|C_i} \cong \omega_{C_i}(\de)$, where
$\de$ is the scheme theoretic intersection of $C_i$ and $\bigcup_{i \neq j}C_j$. It follows that $\omega_{C_i}
\cong \oc_{C_i}(-a+1-\de)$ and since $deg(\de) \geq 0$, this implies
that $p_a(C_i) < 0$, then $C_i$ is a multiple structure on
$\tilde{C_i}$ of multiplicity $>1$.\\
It turns out that each irreducible component of $Z_{red}$ appears with multiplicity $>1$, thus since $deg(Z)=4$ it follows that $Z$ is a double structure on a quadric surface or a 4-uple structure on a plane. This last case can be readily solved. Indeed $C$ would be a 4-uple structure on a line and thanks to \cite{BF} (Remark 4.4) we know that a thick and l.c.i. 4-uple structure on a line is a global complete intersection. Hence we can assume $Z$ quasi-primitive, i.e. we can assume $Z$ does not contain the first infinitesimal neighbourhood  of $Z_{red}$.
Anyway by \cite{M} (see main theorem and Section B) and since $Z_{red}$ is a plane we also have that $Z$ is a c.i..\\
We then suppose that $Z$ is a double structure on a quadric
surface of rank $\geq 2$, which is a complete intersection
$(1,2)$. By Lemma \ref{lem-Z} it follows that $Z$ is a c.i..\epf
\begin{definition}
Let $\ec$ be a rank two normalized vector bundle (i.e. $c_1(\ec)=-1, 0$), we set $r:=min\{n | h^0(\ec(n))\neq 0\}$. If $r>0$, $\ec$ is stable. If $r \leq 0$ we call $r$ \textit{degree of instability of $\ec$}.
\end{definition}
\begin{remark}
The next lemma represents a slight improvement of previous results about the existence of rank two vector bundles in $\Pq$ and $\Pcq$.\\
Indeed Decker proved that any stable rank two vector bundle on $\Pq$ with $c_1=-1$ and $c_2=4$ is isomorphic to the Horrocks-Mumford bundle and that in $\Pcq$ there is no stable rank two vector bundle with these Chern classes (see \cite{D}). We show that neither are there such vector bundles with $r=0$. As for bundles with $c_1=0$ and $c_2=3$, there are similar results by Barth-Elencwajg (see \cite{BE}) and Ballico-Chiantini (see \cite{BC}) stating that $r < 0$. We prove that in fact $r < -1$.
\end{remark}
\begin{lemma}
\label{lem-E}
There does not exist any rank two vector bundle $\ec$ on $\Pq$ such that $r=0$, $c_1(\ec)=-1$, $c_2(\ec)=4$ or, respectively, $r=-1$, $c_1(\ec)=0$, $c_2(\ec)=3$.
\end{lemma}
\textit{Proof:} We observe first of all that in both cases there are no integers $a, b$ satisfying the equations $a+b=c_1$, $ab=c_2$, hence the vector bundle $\ec$ cannot be split.\\
Assume $\ec$ has $r=0$, $c_1(\ec)=-1$, $c_2(\ec)=4$, then $h^0(\ec) \neq 0$. There is a section of $\ec$ vanishing on a codimension two scheme $Z$: $0 \to \oc \to \ec \to \id_{Z}(-1) \to 0$. We have $deg(Z)=c_2(\ec)=4$ and $Z$ subcanonical with $\omega_Z \cong \oc_Z(-6)$.\\
If $r=-1$, $c_1(\ec)=0$, $c_2(\ec)=3$, then $h^0(\ec(-1)) \neq 0$ and we get a section of $\ec(-1)$ vanshing in codimension two along a quartic surface $Z$, with $\omega_Z \cong \oc_Z(-7)$.\\
It is enough to apply \ref{lem-Z2} to conclude that such vector
bundles cannot exist.\epf


\section{Codimension two subvarieties in $\Pn$, $n \geq 5$}
Let $X \subset \Pn$, $n \geq 5$ be a smooth codimension two subvariety, lying on a hypersurface $\sig$ of degree $m \geq 5$ with a ($m$-2)-uple linear subspace $K$ of codimension two, i.e. $K \cong \mathbb{P}^{n-2}$. If $n=5$ we assume $Pic(X)=\bZ H$, for $n \geq 6$ this is granted by Barth's theorem. In any case we set $\omega_X \cong \oc_X(e)$.\\
The general $\Pq$ section $S$ of $X$ is a surface lying on a threefold $\sig \cap H$ of degree $m$ having a singular plane of multiplicity ($m$-2). We will always suppose that $h^0(\id_S(2))=0$.\\
We will prove that $S$ contains a plane curve. First we fix some notations and state some results concerning surfaces containing a plane curve, proofs and more details can be found in \cite{ElFo}.\\
Let $P$ be a plane curve of degree $p$, lying on a smooth surface $S \subset \Pq$. Let $\Pi$ be the plane containing $P$ and let $Z:=S \cap \Pi$. We assume that $P$ is the one-dimensional part of $Z$ and we define $\rc$ as the residual scheme of $Z$ with respect to $P$, namely $\id_{\rc}:=(\id_Z : \id_P)$.
The points of the zero-dimensional scheme $\rc$ can be isolated as well as embedded in $P$.\\
Let $\delta$ be the $\infty^1$ linear system cut out on $S$, residually to $P$, by the hyperplanes containing $\Pi$.
Severi's theorem states that unless $S$ is a Veronese surface, then $h^1(\id_S(1))=0$ and thus $H^0(\oc_{\Pq}(1)) \cong H^0(\oc_S(1))$. Moreover if $p \geq 2$, the hyperplanes containing $\Pi$ are exactly those containing $P$.
This allows us to conclude that $\delta=|H-P|$ (on $S$). We will denote by $Y_H$ the element of $\delta$ corresponding to the hyperplane $H$ and we call $C_H = P \cup Y_H=S \cap H$.\\
Let $\bc$ be the base locus of $\delta$. We have the following results.
\begin{lemma}
\label{lem-P}
(i) P is reduced, the base locus $\bc$ is contained in $\Pi$ and $dim(\bc) \leq 0$. The general $Y_H \in \delta$ is smooth out of $\Pi$ and does not have any component in $\Pi$.\\
(ii) $\bc = \rc$ and $\deg(\rc)=(H-P)^2=d-2p+P^2$.
\end{lemma}
\textit{Proof:} See Lemma 2.1 and 2.4 of \cite{ElFo}.\epf\\
\\
In the present situation, $S$ is subcanonical with $\omega_S \cong \oc_S(e+n-4)$.
We know $deg(\rc)=d-2p+P^2$ and we compute $P^2$ by adjunction, knowing $p_a(P)$ since $P$ is a plane curve and recalling that $K_S=(e+n-4)H$.
It turns out that $deg(\rc)=d+p^2-p(e+n+1)$.
\begin{lemma}
\label{lem-exist}
If $S \subset \Pq$ is a smooth surface, lying on a degree $m$ hypersurface $\sig$ with a $(m-2)$-uple plane, then $S$ contains a (reduced) plane curve, $P$. If $H$ is a general hyperplane through $P$, then $H \cap S=P \cup Y_H$ where $Y_H$ has no irreducible components in $\Pi$ and is smooth out of $\Pi$.
\end{lemma}
\textit{Proof:} If $\Pi$ is the plane with multiplicity $(m-2)$ in $\sig$ and $H$ is an hyperplane containing $\Pi$, we have $H \cap \sig=(m-2)\Pi \cup Q_H$, where $Q_H$ is a quadric surface and $C_H=S \cap H \subset (m-2)\Pi \cup Q_H$. If $dim(C_H \cap \Pi)=0$, then $C_H \subset Q_H$, but this is excluded by our assumptions. Indeed by Severi's theorem $h^0(\id_{C_H}(2)) \neq 0$ would imply $h^0(\id_S(2)) \neq 0$. So $dim(C_H \cap \Pi)=1$ and $S$ contains a plane curve. We conclude with Lemma \ref{lem-P}. \epf\\
\\
If $H$ is an hyperplane through $\Pi$, the corresponding section is $C_H = Y_H \cup P$. Since $Y_H$ does not have any component in $\Pi$, we have $Y_H \subset Q_H$. We denote by $q_H$ the conic $Q_H \cap \Pi$. As $H$ varies, the $q_H$'s form a family of conics in $\Pi$. Let $\bc_q$ be the base locus of $\{q_H\}$, we have $\rc \subset \bc_q$, since $Y_H \cap \Pi \subset Q_H \cap \Pi = q_H$.\\
One can show that $\bc_q$ is $(m-1)$-uple in $\sig$ (see \cite{ElFo}, Lemma 3.3). To prove this, just consider an equation $\varphi$ of $\sig$ and note that clearly $\varphi \in \mathbb{I}^2(\Pi)$. Easy computations show that all $(s-2)$-th derivatives of $\varphi$ vanish at a point $x \in \bc_q$.\\
The following result concerns in particular subcanonical surfaces.
\begin{lemma}
\label{lem-sub}
With notations as above ($S$ subcanonical with $\omega_S \cong \oc_S(a)$), we have:\\
(i) $deg(P) \leq a+3$.\\
(ii) If $\rc=\emptyset$, then $S$ is a complete intersection.
\end{lemma}
\textit{Proof:} (i) We have already computed $deg(\rc)=-p(a+5)+d+p^2$. Recall that $deg(Y_H)=d-p$ and $deg(\rc) \leq deg(Y_H)$, this implies $p \leq a+4$. We will see that the case $p=a+4$ is not possible. Let $p=a+4$, then $Y_H \cdot P=p-P^2=-p(p-a-4)=0$, i.e. $Y_H \cap P=\emptyset$. In other words the curve $C_H=S \cap H = Y_H \cup P$ is not connected, but this is impossible since $h^0(\oc_{C_H})=1$ (use $0 \to \oc_S(-1) \to \oc_S \to \oc_{C_H} \to 0$ and $h^1(\oc_S(-1))=h^1(\omega_S(1))=0$ by Kodaira).\\
(ii) Since $S$ is subcanonical we can consider the exact sequence $0 \to \oc  \to \ec \to
\id_S(a+5) \to 0$. If we restrict it to $\Pi$ and divide by an
equation of $P$, we get $0 \to \oc_{\Pi} {\to} \ec_{\Pi}(-p) \to
\id_{\rc}(a+5-2p) \to 0$. If $\rc=\emptyset$, then
$\id_{\rc}=\oc_{\Pi}$ and the above sequence splits. It follows
that $\ec$ splits and $S$ is a complete intersection.\epf
\begin{example}
Let $S$ be a smooth section of the Horrocks-Mumford bundle $\fc$, $S$ is an abelian variety and has $\omega_S=\oc_S$. By Lemma \ref{lem-sub} we know that if $S$ contains a plane curve $P$, then $p \leq 3$. Moreover, $P$ cannot be a line or a conic, since these curves are rational and this would imply that there exists a non constant morphism $\Pun \to S$, factoring through $Jac^0(\Pun) \cong \{\ast\}$ and this is not possible. Then necessarily $P$ is a plane smooth cubic (hence elliptic).\\
By the "reducibility lemma" of Poincar\'e, an abelian surface $S$ contains an elliptic curve if and only if $S$ is isogenous to a product of elliptic curves. It is known that the general section of the Horrocks-Mumford bundle is not isogenous to a product of elliptic curves, but there exist smooth sections satisfying such property (see \cite{R}, \cite{L}). Summarizing we can say that among the sections of Horrocks-Mumford bundle we can find smooth surfaces containing a plane curve, but the general one does not contain any.\\
Now assume $S$ to be one of those smooth surfaces containing a plane cubic, $P$. Let $\Pi$ be the plane spanned by $P$. Recall that we have $0 \to \oc \stackrel{s}{\to} \fc(3) \to \id_S(5) \to 0$. We restrict the sequence to $\Pi$ and since $s_{|\Pi}$ vanishes along $P$, we can divide by an equation of $P$ and obtain a section of $\fc_{|\Pi}$. We then have $h^0(\fc_{|\Pi}) \neq 0$, i.e. $\fc_{|\Pi}$ is not stable, in other words $\Pi$ is an unstable plane for $\fc$.
\end{example}
In order to prove Theorem \ref{th} we need some other preliminary results.
\begin{lemma}
\label{lem-min-deg}
Let $F \subset \Pt$ be a surface of degree $m$, singular along a line $D$ with multiplicity $m-1$. Then $F$ is the projection of a surface of degree $m$ in $\mathbb{P}^{m+1}$ (minimal degree surface).
\end{lemma}
\textit{Proof:} The surface $F$ is rational. Let $p: F' \to F$ be a desingularization of $F$ and let $H$ be a divisor in $p^*\oc_F(1)$. We have $0 \to \oc_{F'} \to \oc_{F'}(H) \to \oc_H(H) \to 0$ and since $F'$ is rational too, then $h^1(\oc_{F'})=0$. Now $h^0(\oc_H(H))=m+1$ ($H$ is a rational curve), then $h^0(\oc_{F'}(H))=m+2$ and we can embed $F'$ in $\mathbb{P}^{m+1}$.\epf
\begin{remark}
Minimal degree surfaces in $\Pn$ are classified, in particular they can be: a smooth rational scroll, a cone over a rational normal curve of $\mathbb{P}^{n-1}$ or the Veronese surface if $n=5$. Except for the Veronese, all these surfaces are ruled in lines.
\end{remark}
\begin{lemma}
\label{lem-ruled}
Let $T \subset \mathbb{P}^{m+1}$, $m \geq 3$, be a surface ruled in lines. Let $C \subset T$ be a smooth irreducible curve. If $dim(<C>)=3$, then $deg(C) \leq deg(T)-m+3$. ($<C>$ is the linear space spanned by $C$)
\end{lemma}
\textit{Proof:} Let us consider $m-3$ general points on $C$ and let $f_1,\ldots,f_{m-3}$ be the rulings passing through these points. We consider moreover $m-3$ points $p_1,\ldots,p_{m-3}$ such that $p_i \in f_i$ but $p_i \not \in <C>$ and let also $q_1,\ldots,q_4$ be four general points in $<C>$. We thus have $m+1$ points, spanning at most a space of dimension $m$, hence these points are contained in a hyperplane $H$ of $\mathbb{P}^{m+1}$. Now $<C> \subset H$ since $q_i \in H$ $\forall$ $i=1,\ldots,4$, $f_i \subset H$ since $card(f_i \cap H)>1$ $\forall$ $i=1,\ldots,m-3$, so $H \cap T$ contains $C,f_1,\ldots,f_{m-3}$ (which form a degenerate curve in $T$ of degree $m-3+deg(C)$) and this yields: $deg(T) \geq deg(C)+m-3$.\epf\\
\begin{lemma}
\label{lem-hyper}
Let $X,K \subset \Pn$, $n \geq 4$, $X$ smooth of codimension two, $K \cong \mathbb{P}^{n-2}$ a linear subspace. Let $dim(X \cap K)=n-3$. If the general hyperplane section of $X \cap K$ contains a linear subspace of dimension $n-4$, then $X$ contains a linear subspace of dimension $n-3$.
\end{lemma}
\textit{Proof:} We see $X_K=X \cap K$ as a hypersurface in $K \cong \mathbb{P}^{n-2}$. A general hyperplane of $K$ is cut on $K$ by a general hyperplane of $\Pn$. Then the hypersurface $X_K$ of $K$ is such that its general hyperplane section contains a linear subspace of dimension $n-4$. We claim that $X_K$ contains an hyperplane of $K$. Indeed we may assume $X_K$ reduced. Let $X_K=T_1 \cup \ldots \cup T_r$ be the decomposition of $X_K$ into irreducible components. Now using the fact that the general hyperplane section of each $T_i$ is irreducible, we conclude that one of the $T_i$'s has degree one and thus $X_K$ contains an hyperplane of $K$.\epf\\
\\
\textit{Proof of Theorem \ref{th}:} We only need to work out the case $n=5$. We will follow the method used in the proof of Theorems 1.1 and 1.2 of \cite{ElFo}. We must distinguish three cases, depending on the behaviour of the base locus $\bc_q$ of the conics $q_H$. If $dim(\bc_q)=0$, at least two of the conics intersect properly and then $deg(\bc_q) \leq 4$. It follows that $r:=deg(\rc) \leq 4$ too, since $\rc \subset \bc_q$. If $dim(\bc_q)=1$, there are two possibilities: the one-dimensional part of $\bc_q$ can be a line or a conic.\\
If the conics $q_H$ move, i.e. if $dim(\bc_q)=0$, we have seen that $r = deg(\rc) \leq 4$. We observe that actually we can suppose $r \geq 1$, indeed by \ref{lem-sub} if $\rc= \emptyset$, then $S$ (and $X$) is a complete intersection.\\
If $H$ is a general hyperplane, $Y_H \cap P \subset q_H \cap P$ and since at least one conic intersects $P$ properly, we obtain $Y_H.P \leq 2p$. We have $Y_H.P=p-P^2$ and recalling that $r=d-2p+P^2$, it follows $Y_H.P=d-p-r$. Putting everything together: $p \geq \frac{d-r}{3} \geq \frac{d-4}{3}$. On the other hand we have $Y_H.P=p(e+5-p)$ and clearly this implies $p \geq e+3$. Comparing this with the result stated in \ref{lem-sub} and setting $\omega_S \cong \oc_S(e+1)$, we are left with only two possibilities: $p=e+3$ or $p=e+4$. We have already observed that $d=p(e+6)-p^2+r$, then considering the two cases above, we can express $d$ in terms of $e$ and $r$ and we get the following formulas:
$$if \:\: p=e+3,\: then\:\: d=3(e+3)+r \:\:\:\:\:\:\:\:\:\:(2)$$
$$if \:\: p=e+4,\: then\:\: d=2(e+4)+r \:\:\:\:\:\:\:\:\:\:(3)$$
We recall that if $C$ lies on a quartic surface and $d$ is large enough, $X$ lies on a quartic hypersurface too, then $X$ is a complete intersection. We know that $\pi-1=\frac{d(e+2)}{2}$, then since $ \frac{d-4}{3} \leq p \leq e+4$ we obtain $\pi-1 \geq \frac{d(d-10)}{6}$. If we compare this quantity with $G(d,5)$, we see that if $d \geq 33$, then $h^0(\id_C(4)) \neq 0$ and $X$ is a complete intersection.\\
Thanks to the result in Proposition \ref{prop} we know that if $d \leq 25$, $X$ is a complete intersection too, then we only have to check the cases $26 \leq d \leq 32$.\\
We assume $h^0(\id_C(4))=0$, then it must be $\pi=1+\frac{d(e+2)}{2} \leq G(d,5)$. Thanks to this inequality it is easy to see that for $d \leq 32$, we always have $e \leq 5$. Now if we look at formulas $(2)$ and $(3)$ above, clearly $e \leq 5$ implies $d \leq 28$.\\
On the other hand, in order to have $d \geq 26$, $e$ must be at least equal to $4$.\\
If $d=26,27,28$, the condition on the genus $\pi$ yields $e=4$ again. However, if we look at formulas $(2)$ and $(3)$ we see that if $e=4$, $d$ is at most equal to $25$.\\
If $\bc_q$ contains a line, $D$, then $D$ has multiplicity $m-1$ in $\sig$, so if $H$ is an hyperplane containing $D$ (but not $\Pi$), $F=\sig \cap H$ is a surface of degree $m$ in $\Pt$ having a $(m-1)$-uple line. This kind of surface is a projection of a degree $m$ surface in $\mathbb{P}^{m+1}$, by Lemma \ref{lem-min-deg}. The hyperplane section $C=S \cap H$ is a curve contained in $F$. We must distinguish two cases: $D \subset S$ or $D \not \subset S$.\\
If $D \not \subset S$, we claim that the general $C$ is smooth. Let $|L|$ be the linear system cut on $S$ by the hyperplanes containing $D$ and let $B=D \cap S=\{p_1,\ldots,p_r\}$. Since $B$ is the base locus of $|L|$, the general element of $|L|$ is smooth out of $B$. If all curves in $|L|$ were singular at a point $p_i \in B$, it would be $T_{p_i}S \subset H$, $\forall H \supset D$. Anyway the intersection of $H \supset D$ is only $D$, so this is not possible.
It follows that the curves of $|L|$ singular at a $p_i \in B$ form a closed subset of $|L|$. The same holds for all $p \in B$, hence the claim.\\
Let $F'$ be a surface in $\mathbb{P}^{m+1}$ projecting down to $F$. Since $C$ is not contained in the singular locus of $F$, there exists a curve $C' \subset F'$ such that the projection restricted to $C'$ is an isomorphism over $C$. In particular $\oc_{C'}(1) \cong \oc_C(1)$ and since $C$ is linearly normal in $\Pt$, this implies that $C'$ is degenerate. Now we can apply Lemma \ref{lem-ruled} to $F'$ and $C'$ (we have already pointed out that $F'$ is ruled in lines unless $F'$ is the Veronese surface) and we get $d=deg(C') \leq m-m+3=3$. If $F'$ is the Veronese surface we have anyway $d \leq 4$.\\
If $D \subset S$, then $D$ is a component of the plane curve $P$ (the one-dimensional part of $S \cap \Pi$). Coming back to the variety $X \subset \sig \subset \Pcq$ with $K \cong \Pt \subset \sig$ a linear subspace of multiplicity $m-2$, we have a surface $X_K=X \cap K \subset K \cong \Pt$ such that its general hyperplane section contains a line. This implies by Lemma \ref{lem-hyper} that $X_K$ contains a plane and thus $X$ contains a plane, say $E$. This plane is a Cartier divisor on the smooth threefold $X$. Since we are supposing $Pic(X) = \bZ H$, there exists an hypersurface such that $E$ is cut on $X$ by this hypersurface, but this could happen only if $deg(X)=1$.\\
To complete the proof we only have to consider the case in which $\bc_q=q$, where $q$ is an irreducible conic (if $q$ is reducible, $\bc_q$ contains a line). For every $Y_H \in |H-P|$ we consider the zero-dimensional scheme $\de_H=Y_H \cap q$. For every $H$, $\de_H$ is a subset of $d-p$ points of $q$.\\
There are two possibilities: $q \subset S$ or $q \not \subset S$. If $q \not \subset S$, then $\de_H$ is fixed (otherwise the points of $\de_H$ would cover the conic, as $H$ varies, i.e. $q \subset S$). It must be $\de_H= \rc$. It is enough to compare the degrees of $\de_H$ and $\rc$ to see that this implies $P^2=p$ and then $Y_H.P=P^2-p=0$. This is not possible since the corresponding hyperplane section $C_H$ of $S$ would be disconnected.\\
Hence $q \subset S$ and then $q \subset P$. In other words: $\de_H=Y_H \cap P$, thus $Y_HP=d-p$ and $r=0$. By Lemma \ref{lem-sub} we conclude that $X$ is a complete intersection.\epf\\

\bigskip
\noindent
Address of the author:
\par \noindent
Dipartimento di Matematica
\par \noindent
35, via Machiavelli
\par \noindent
44100 Ferrara (Italy)
\par \noindent
email:  chiara@dm.unife.it\\

\end{document}